\documentclass{amsart}

\usepackage{xypic}
\input xy
\xyoption{all}
\usepackage{epsfig}
\usepackage{amsthm}
\usepackage{amssymb}
\usepackage{amsmath}
\usepackage{amscd}
\usepackage{color}
\usepackage{mathabx}
\usepackage[T1]{fontenc}

\newcommand{\bg}{\begin{equation}}
\newcommand{\ed}{\end{equation}}
\newcommand{\bga}{\begin{eqnarray}}
\newcommand{\eda}{\end{eqnarray}}
\newcommand{\pf}{\textbf{Proof:\ }}

\def\cbdu{\par{\raggedleft$\Box$\par}}

\newtheorem {Theorem}  {Theorem}

\numberwithin{Theorem}{section}

\newtheorem {Lemma}[Theorem]  {Lemma}

\theoremstyle{definition}

\theoremstyle{remark}
\newtheorem{Remark}[Theorem]{\bf Remark}

\newtheorem {Corollary}[Theorem]{\bf Corollary}

\expandafter\chardef\csname pre amssym.def
at\endcsname=\the\catcode`\@ \catcode`\@=11
\def\undefine#1{\let#1\undefined}
\def\newsymbol#1#2#3#4#5{\let\next@\relax
 \ifnum#2=\@ne\let\next@\msafam@\else
 \ifnum#2=\tw@\let\next@\msbfam@\fi\fi
 \mathchardef#1="#3\next@#4#5}
\def\mathhexbox@#1#2#3{\relax
 \ifmmode\mathpalette{}{\m@th\mathchar"#1#2#3}%
 \else\leavevmode\hbox{$\m@th\mathchar"#1#2#3$}\fi}
\def\hexnumber@#1{\ifcase#1 0\or 1\or 2\or 3\or 4\or 5\or 6\or 7\or 8\or
 9\or A\or B\or C\or D\or E\or F\fi}

\font\teneufm=eufm10 \font\seveneufm=eufm7 \font\fiveeufm=eufm5
\newfam\eufmfam
\textfont\eufmfam=\teneufm \scriptfont\eufmfam=\seveneufm
\scriptscriptfont\eufmfam=\fiveeufm

\catcode`\@=\csname pre amssym.def at\endcsname

\newcounter{remark}
\newcommand{\R}{\mathbf{R}}

\def  \R   {{\mathbb R}}

\def  \12  {{\frac{1}{2}}}

\def\build#1_#2^#3{\mathrel{\mathop{\kern 0pt#1}\limits_{#2}^{#3}}}

\begin{document}

\title[Decay of Hall-MHD]{Long time behavior of solutions to the 3D Hall-magneto-hydrodynamics system with one diffusion}

\author [Mimi Dai]{Mimi Dai}
\address{Department of Mathematics, Stat. and Comp.Sci., University of Illinois Chicago, Chicago, IL 60607,USA}
\email{mdai@uic.edu} 

\author [Han Liu]{Han Liu}
\address{Department of Mathematics, Stat. and Comp.Sci., University of Illinois Chicago, Chicago, IL 60607,USA}
\email{hliu94@uic.edu}

\thanks{The authors were partially supported by NSF grant
DMS--1815069.}





\begin{abstract}
This paper studies the asymptotic behavior of smooth solutions to the generalized Hall-magneto-hydrodynamics system (\ref{HMHD}) with one single diffusion on the whole space $\R^3$. We establish that, in the inviscid resistive case, the energy $\|b(t)\|_2^2$ vanishes and $\|u(t)\|_2^2$ converges to a constant as time tends to infinity provided the velocity is bounded in $W^{1-\alpha,\frac3\alpha}(\R^3)$; in the viscous non-resistive case, the energy $\|u(t)\|_2^2$ vanishes and $\|b(t)\|_2^2$ converges to a constant provided the magnetic field is bounded in $W^{1-\beta,\infty}(\R^3)$. In summary, one single diffusion, being as weak as $(-\Delta)^\alpha b$ or $(-\Delta)^\beta u$ with small enough $\alpha, \beta$, is sufficient to prevent asymptotic energy oscillations for certain smooth solutions to the system.

\bigskip


\hspace{0.02cm}CLASSIFICATION CODE: 35Q35, 35B40, 35Q85.
\end{abstract}

\maketitle

\section{Introduction}

In this paper we consider the three dimensional incompressible 
Hall-magneto-hydrodynamics (Hall-MHD) system with fractional Laplacian:
\begin{equation}\label{HMHD}
\begin{split}
u_t+u\cdot\nabla u-b\cdot\nabla b+\nabla p=-\nu(-\Delta)^\beta u,\\
b_t+u\cdot\nabla b-b\cdot\nabla u+\eta\nabla\times((\nabla\times b)\times b)=-\mu(-\Delta)^\alpha b,\\
\nabla \cdot u=0, 
\end{split}
\end{equation}
with the initial conditions
\begin{equation}
u(x,0)=u_0(x),\qquad b(x,0)=b_0(x), \qquad \nabla\cdot u_0=\nabla\cdot b_0=0,
\end{equation}
where $x\in\mathbb{R}^3$, $t\geq 0$; the unknown functions $u, p, b$ denote the fluid velocity, fluid pressure and magnetic field, respectively. The constants $\nu, \mu$ stand for the kinematic viscosity coefficient of the fluid and the magnetic Reynolds number, respectively. The Hall parameter $\eta$ denotes the strength of magnetic reconnection, which takes value either $1$ or $0$ in this paper. We point out that the magnetic field $b$ will remain divergence free for all the time if $\nabla\cdot b_0=0$, see \cite{CL}. 
When the parameters $\alpha=\beta=\eta=1$, the Hall term $\nabla\times((\nabla\times b)\times b)$ describes the occurrence of the magnetic reconnection when the magnetic shear is large, which makes (\ref{HMHD}) different from the usual MHD system. The magnetic reconnection is a physical process of topological reordering of magnetic field lines, in which magnetic energy is transferred to kinetic energy, thermal energy and particle acceleration.
In this paper, we consider the generalized Hall-MHD system (\ref{HMHD}) with $0<\alpha,\beta\leq 1$ and $\eta=1$. The readers are referred to \cite{For, Li, SC} and references therein for more physical background of the Hall-MHD system.

We will go over briefly the mathematical study on the Hall-MHD system ($\alpha=\beta=\eta=1$) and some fundamental theory in the literature. 
Global existence of weak solutions was established in \cite{ADFL} and \cite{CDL,DS}, on periodic domain and in the whole space $\R^3$ respectively.  Local well-posedness of classical solution was investigated in \cite{CDL}. Various blow-up criteria were obtained in \cite{CDL,CL,Dai}. Local well-posedness is studied in \cite{CWW} for the Hall-MHD system with fractional magnetic diffusion. In \cite{CS} the authors established temporal decay estimates for weak solutions in the energy space $L^2(\R^3)$ and for small initial data solutions in higher order Sobolev spaces by utilizing the Fourier splitting technique introduced in \cite{S0,S1,S3,S2}, and extensively applied to study asymptotic behavior of solutions to various systems on the whole space $\R^3$, for instance see \cite{DFRSS, DQS, DSch}. Specifically, provided the initial data $(u_0,b_0)\in\left(L^2(\R^3)\cap L^1(\R^3)\right)^2$, there exists a weak solution to system (\ref{HMHD}) with $\alpha=\beta=1$ which satisfies 
\begin{equation}\label{op-decay}
\|u(t)\|_2^2+\|b(t)\|_2^2\leq C(1+t)^{-\frac32}.
\end{equation}
Notice that for the Navier-Stokes equation, the energy $\|u(t)\|_2^2$ of a weak solution decays with the same rate $(1+t)^{-\frac32}$ if the initial data $u_0\in L^2(\R^3)\cap L^1(\R^3)$, see \cite{S1}. 
In contrast, in \cite{AS} the authors studied the long time behavior of solutions to the MHD system, that is (\ref{HMHD}) with $\alpha=\beta=1$ and $\eta=0$ (thus without the Hall term $\nabla\times((\nabla\times b)\times b)$). In the viscous and resistive case $\mu,\nu>0$, it was shown that any weak solution $(u,b)$ to the MHD system with initial data in $\left(L^2(\R^3)\right)^2$, the total energy $\|u(t)\|_2^2+\|b(t)\|_2^2$ converges to zero without a rate as time tends to infinity. 
The authors also analyzed that this decay estimate without rate is actually optimal. One can expect that if we assume additionally that $(u_0,b_0)$ is also in $\left(L^1(\R^3)\right)^2$, then the energy $\|u(t)\|_2^2+\|b(t)\|_2^2$ would decay with a rate $(1+t)^{-\frac32}$ which is the same as (\ref{op-decay}). It indicates that the presence of the Hall term in (\ref{HMHD}) seems to make no difference for the energy decay estimates when we compare the viscous resistive Hall-MHD and MHD systems. Another case considered in \cite{AS} was when $\nu>0,\mu=0$, it was shown that if strong bounded solutions to the viscous non-resistive MHD system exist, then the diffusion of the fluid velocity is enough to prevent any asymptotic energy oscillations, thus $\|u(t)\|_2$ converges to zero and $\|b(t)\|_2$ converges to a constant as $t\to\infty$.

Inspired by the work of \cite{AS, CS}, we are interested in studying the asymptotic behavior of smooth solutions to the generalized Hall-MHD system (\ref{HMHD}) ($\eta=1$) with solely one diffusion, that is, either the viscous non-resistive case $\nu>0,\mu=0$ or the inviscid resistive case $\nu=0,\mu>0$.
The subjects of the investigation are twofold:  to explore the different or analogous large time behavior between solutions of the Hall-MHD system and the MHD system; to examine how strong the one diffusion is needed to prevent asymptotic energy oscillations. To carry out the analysis, we split the energy into high and low frequency parts and estimate them separately.  It turns out that the energy of the low frequency part converges once the initial data has finite energy; while additional assumptions are necessary to enforce the convergence of the energy of the high frequency part. The Fourier splitting approach is applied to establish the decay of the high frequency part. The main results are stated as below.


\begin{Theorem}\label{thm1}
Let $(u,b)$ be a global smooth solution to (\ref{HMHD}) with $\nu =0$ and $\mu>0$. Assume $u_0\in L^2(\R^3)$ and $b_0\in L^1(\R^3)\cap L^2(\R^3)$, and one of the two conditions holds: 
\begin{equation}\notag
\begin{split}
(i)\ \ &  u\in L^\infty(0,\infty;W^{1-\alpha, \frac3\alpha}(\R^3)) \text{ with }  \alpha\in[1/2, 1]; \\
(ii)\ \ & b \in L^\infty(0,\infty;W^{1-\alpha,\infty}(\R^3)) \text{ and } u\in L^\infty(0,\infty;W^{1-\alpha, \frac3\alpha}(\R^3)) \text{ with } \alpha\in (0,1].
\end{split}
\end{equation}
 Then we have 
 \[\lim_{t\to\infty} \|b(t)\|^2_2 =0, \ \ \lim_{t\to\infty}\|u(t)\|^2_2 =C\]
 for some absolute constant $C$.
\end{Theorem}

\begin{Theorem}\label{thm2}
Let $(u,b)$ be a global smooth solution to (\ref{HMHD}) with $\mu =0$, $\nu>0$ and $0<\beta\leq 1$. Assume $u_0\in L^1(\R^3)\cap L^2(\R^3)$ and $b_0\in L^2(\R^3)$, and additionally $b \in L^\infty(0,\infty;W^{1-\beta,\infty}(\R^3))$.
 Then we have 
 \[\lim_{t\to\infty} \|u(t)\|^2_2 =0, \ \ \lim_{t\to\infty}\|b(t)\|^2_2 =C\]
 for some absolute constant $C$.
\end{Theorem}

As stated in Theorem \ref{thm1} for the inviscid resistive Hall-MHD system, we discovered that if the velocity is bounded in $W^{1-\alpha, \frac3\alpha}(\R^3)$, then the single diffusion $(-\Delta)^{\alpha}b$ with $\alpha\geq \frac12$ is enough to guarantee the convergence of the two energies and thus prevent asymptotic energy oscillations. Furthermore, if both the velocity and the magnetic field are bounded in some higher order Sobolev space, then the single diffusion $(-\Delta)^{\alpha}b$ with $\alpha>0$ is sufficient to prevent the asymptotic energy oscillations. As a contrast, for the viscous non-resistive Hall-MHD system in Theorem \ref{thm2}, if the magnetic field $b$ is bounded in $W^{1-\beta,\infty}(\R^3)$, then the single diffusion $(-\Delta)^{\beta}u$ with $\beta>0$ is adequate to stave off such oscillations. We point out that the additional assumptions are imposed to estimate the nonlinear terms in which no cancelation exists.

In addition, we learned that no extra assumption is needed to estimate the Hall-term in both cases, thanks to the cancelation in the flux $\int_{\R^3} \nabla\times((\nabla\times b)\times b)\cdot b\, \mathrm dx$.
Therefore, we have the following analogous results for the MHD system.

\begin{Corollary}\label{cor1}
Let $(u,b)$ be a global smooth solution to the generalized inviscid resistive MHD system (\ref{HMHD}) with $\eta=0$, $\nu =0$ and $\mu>0$. Assume $u_0\in L^2(\R^3)$ and $b_0\in L^1(\R^3)\cap L^2(\R^3)$, and either (i) or (ii) in Theorem \ref{thm1} holds.
 Then we have 
 \[\lim_{t\to\infty} \|b(t)\|^2_2 =0, \ \ \lim_{t\to\infty}\|u(t)\|^2_2 =C\]
 for some absolute constant $C$.
\end{Corollary}

\begin{Corollary}\label{cor2}
Let $(u,b)$ be a global smooth solution to the generalized viscous non-resistive MHD system (\ref{HMHD}) with $\eta=0$, $\mu =0$, $\nu>0$ and $0<\beta\leq 1$. Assume $u_0\in L^1(\R^3)\cap L^2(\R^3)$ and $b_0\in L^2(\R^3)$, and additionally $b \in L^\infty(0,\infty;W^{1-\beta,\infty}(\R^3))$.
 Then we have 
 \[\lim_{t\to\infty} \|u(t)\|^2_2 =0, \ \ \lim_{t\to\infty}\|b(t)\|^2_2 =C\]
 for some absolute constant $C$.
\end{Corollary}

\begin{Remark} {\color{white}.}
\begin{itemize}
\item
Corollary \ref{cor2} recovers the result in \cite{AS} when $\beta=1$.
\item
The cancelation of nonlinear terms, as $\langle u\cdot\nabla u,u\rangle=0$, $\langle u\cdot\nabla b,b\rangle=0$, and $\langle \nabla \times ((\nabla \times b)\times b),b\rangle=0$, play a crucial rule in the analysis. These cancelations are valid for smooth solutions, which is the reason we work with smooth solutions. In fact, it was shown in \cite{DS} that these cancelations are also valid for solutions in ``Onsager" space, that is, $u\in L^3(0,\infty;\dot B^{\frac13}_{3, c_0}(\R^3))$ and $b\in L^3(0,\infty;\dot B^{\frac23}_{3, c_0}(\R^3))$. Thus, the results in the theorems and corollaries above hold for solutions $(u,b)\in L^3(0,\infty;\dot B^{\frac13}_{3, c_0}(\R^3))\times L^3(0,\infty;\dot B^{\frac23}_{3, c_0}(\R^3))$ as well.
\end{itemize}
\end{Remark}

The rest of the paper is organized as follows: in Section \ref{sec:pre} we  introduce some notations and establish certain generalized energy inequalities for low and high frequency parts, as well as some auxiliary estimates to handle the high frequency part; Section \ref{sec:dec} and Section \ref{sec:dec1} are devoted to proving Theorem \ref{thm1} and \ref{thm2}, respectively.  As explained above, the Hall term can be estimated without additional assumptions, in light of which the proof of Corollary \ref{cor1} and Corollary \ref{cor2} will be omitted.

\bigskip

\section{Preliminaries}
\label{sec:pre}

\subsection{Notation}
\label{sec:notation}

We denote by $A\lesssim B$ an estimate of the form $A\leq C B$ with
some absolute constant $C$. We write $\|\cdot\|_p=\|\cdot\|_{L^p(\R^3)}$ for simplification; and $\langle\cdot, \cdot\rangle$ stands for the $L^2$-inner product. 

The Fourier transform of a function $f$ is denoted by 
\[\hat{f}=f^{\wedge}=(2\pi)^{-\frac{n}{2}} \int_{\mathbb{R}^n} e^{-i\xi \cdot x}f(x) \mathrm{d}x,\]
and the inverse Fourier transform of a function $\varphi$ is denoted by
\[\check{\varphi}=(2\pi)^{-\frac{n}{2}} \int_{\mathbb{R}^n} e^{i x \cdot \xi}\varphi(\xi) \mathrm{d}\xi.\]

Various constants shall be denoted by $C$ throughout the paper.

\bigskip

\subsection{Generalized energy inequalities}

We first recall that the following energy equality holds for a regular solution of (\ref{HMHD})  
\begin{equation}\notag
\|u(t)\|_2^2+\|b(t)\|_2^2+2\nu\int_{0}^t\|\nabla^\beta u(s)\|_2^2\,\mathrm ds
+2\mu\int_{0}^t\|\nabla^\alpha b(s)\|_2^2\,\mathrm ds
= \|u(0)\|_2^2+\|b(0)\|_2^2.
\end{equation}

Thanks to the identity
\begin{equation}\notag\label{vec}
(\nabla\times b)\times b=b\cdot\nabla b-\nabla \frac{b\cdot b}2=\nabla \cdot (b \otimes b)-\nabla \frac{|b|^2}2,
\end{equation}
and the fact $\nabla\times (\nabla v)=0$ for any function $v$, the Hall term can be rewritten as
\[\nabla \times((\nabla \times b)\times b)=\nabla\times \nabla\cdot(b\otimes b).\]

In the inviscid resistive case $\nu=0, \mu>0$, we take functions $\varphi(\xi)=e^{-|\xi|^2}$ and $\psi(\xi)=1-\varphi(\xi)$ in the Fourier space.
Obviously, $\varphi \hat b$ and $\psi\hat b$ represent the low and high frequency parts of $\hat b$, respectively. It follows from Plancherel's theorem that 
\[\|b(t)\|_2=\|\hat b(t)\|_2\leq \|\varphi\hat b(t)\|_2+\|\psi \hat b(t)\|_2.\] 
In the viscous non-resistive case $\nu>0, \mu=0$, we choose $\varphi(\xi,t)=e^{-|\xi|^{2\beta}t}$ and $\psi(\xi,t)=1-\varphi(\xi,t)$ instead by abusively using the same letters, and split $\|u(t)\|_2$ as 
\[\|u(t)\|_2=\|\hat u(t)\|_2\leq \|\varphi(t)\hat u(t)\|_2+\|\psi(t) \hat u(t)\|_2.\]
In order to estimate these terms, we first establish certain generalized energy inequalities for the  low and high frequency parts respectively.

\begin{Lemma}
Assume $\mu>0$. Choose function $E(t) \in C^1(\mathbb{R}; \mathbb{R}_{+})$ with $E(t) \geq 0$ and $\varphi(\xi)=e^{-|\xi|^2}$, $\psi(\xi)=1-\varphi(\xi)$. Then a weak solution $(u,b)$ of (\ref{HMHD}) satisfies the generalized energy inequalities,
\begin{equation}\label{elow}
\begin{split}
&\|\widecheck\varphi*b(t)\|_2^2\\
\leq &\|e^{-\mu(t-s)(-\Delta)^\alpha}\widecheck\varphi*b(s)\|_2^2-2\int_s^t\langle u\cdot\nabla b(\tau),e^{-2\mu(\tau-s)(-\Delta)^\alpha}\widecheck\varphi*\widecheck \varphi*b(s)\rangle\,\mathrm d\tau\\
&+2\int_s^t\langle b\cdot\nabla u(\tau),e^{-2\mu(\tau-s)(-\Delta)^\alpha}\widecheck\varphi*\widecheck \varphi*b(s)\rangle\,\mathrm d\tau\\
&-2\int_s^t\langle\nabla\times(\nabla\cdot(b\otimes b)(\tau)), e^{-2\mu(-\Delta)^\alpha (\tau-s)}\widecheck\varphi*\widecheck\varphi* b(s)\rangle\, \mathrm d\tau,
\end{split}
\end{equation}
\begin{equation}\label{ehigh}
\begin{split}
E(t) \| \psi\widehat b(t)\|_2^2\leq& E(s)\|\psi\widehat b(s)\|_2^2-2\mu\int_s^tE(\tau)\||\xi|^\alpha\psi\widehat b(\tau)\|_2^2\,\mathrm d\tau\\
&+\int^t_s E'(\tau) \|\psi\widehat b(\tau)\|_2^2\, \mathrm d\tau-2\int_s^tE(\tau)\langle \widehat{u\cdot\nabla b}(\tau),\psi^2\widehat b(\tau)\rangle\,\mathrm d\tau\\
&+2\int_s^tE(\tau)\langle \widehat{b\cdot\nabla u}(\tau),\psi^2\widehat b(\tau)\rangle\,\mathrm d\tau\\
&-2\int^t_s E(\tau)\langle{\nabla\times(\nabla\cdot(b\times b)(\tau))}^\wedge , \psi^2\widehat b(\tau)\rangle\, \mathrm d\tau.
\end{split}
\end{equation}
\end{Lemma}
\pf
The estimates will be established formally for classical solutions. Multiplying the second equation of \ref{HMHD} by $e^{-2\mu(t-s)(-\Delta)^\alpha}\widecheck\varphi*\widecheck\varphi*b(s)$ and integrating over $\R^3$ yields
\begin{equation}\notag
\begin{split}
&\langle b_t,e^{-2\mu(t-s)(-\Delta)^\alpha}\widecheck\varphi*\widecheck\varphi*b(s)\rangle
+\mu\langle \nabla b, e^{-2\mu(t-s)(-\Delta)^\alpha}\widecheck\varphi*\widecheck\varphi*\nabla b(s)\rangle\\
&+\langle u\cdot\nabla b, e^{-2\mu(t-s)(-\Delta)^\alpha}\widecheck\varphi*\widecheck\varphi*b(s)\rangle
-\langle b\cdot\nabla u, e^{-2\mu(t-s)(-\Delta)^\alpha}\widecheck\varphi*\widecheck\varphi*b(s)\rangle\\
&+\langle \nabla\times(\nabla\cdot (b\times b)),e^{-2\mu(t-s)(-\Delta)^\alpha}\widecheck\varphi*\widecheck\varphi*b(s)\rangle=0.
\end{split}
\end{equation}
We rewrite the first two terms as
\begin{equation}\notag
\begin{split}
&\frac12\frac d{dt}\|e^{-\mu(t-s)(-\Delta)^\alpha}\widecheck\varphi*b(s)\|_2^2-\langle e^{-\mu(t-s)\Delta}\widecheck\varphi*b(s), \partial_t\left(e^{-\mu(t-s)(-\Delta)^\alpha}\widecheck\varphi\right)*b\rangle\\
&+\mu\langle e^{-\mu(t-s)(-\Delta)^\alpha}\widecheck\varphi*b(s), (-\Delta)^\alpha\left(e^{-\mu(t-s)(-\Delta)^\alpha}\widecheck\varphi\right)*b\rangle\\
=&\frac12\frac d{dt}\|e^{-\mu(t-s)(-\Delta)^\alpha}\widecheck\varphi*b(s)\|_2^2
\end{split}
\end{equation}
which is obtained due to the fact $\partial_t\left(e^{-\mu(t-s)(-\Delta)^\alpha}\widecheck\varphi\right)=-\mu(-\Delta)^\alpha \left(e^{-\mu(t-s)(-\Delta)^\alpha}\widecheck\varphi\right)$. Integrating the equation above over the time interval $[s,t]$ produces (\ref{elow}).

We take Fourier transform of the second equation in (\ref{HMHD}), multiply it by $\psi^2\widehat bE(t)$, and integrate over $\R^3$ to infer
\begin{equation}\notag
\begin{split}
&\frac12\frac{d}{dt}\int_{\R^3}E(t)|\psi\widehat b|^2\, \mathrm d\xi-\frac12\int_{\R^3}E'(t)|\psi\widehat b|^2\, \mathrm d\xi+\mu E(t)\int_{\R^3}|\xi|^{2\alpha}|\psi\widehat b|^2\, \mathrm d\xi\\
&+E(t)\langle \widehat {u\cdot\nabla b},\psi^2\widehat b\rangle-E(t)\langle \widehat {b\cdot\nabla u},\psi^2\widehat b\rangle+E(t)\langle  {\nabla\times(\nabla\cdot (b\times b))}^\wedge,\psi^2\widehat b\rangle=0.
\end{split}
\end{equation}
Integrating the last equation over $[s,t]$ yields (\ref{ehigh}).
\cbdu

Similar computation will yield the generalized energy inequalities for the velocity $u$ when $\nu>0$.
\begin{Lemma}
Assume $\nu>0$. Choose function $E(t) \in C^1(\mathbb{R}; \mathbb{R}_{+})$ with $E(t) \geq 0$ and $\varphi(\xi,t)=e^{-|\xi|^{2\beta}t}$, $\psi(\xi, t)=1-\varphi(\xi,t)$. Then a weak solution $(u,b)$ of (\ref{HMHD}) satisfies the generalized energy inequalities,
\begin{equation}\label{u-elow}
\begin{split}
&\|\widecheck\varphi*u(t)\|_2^2\\
\leq &\|e^{-\nu(t-s)(-\Delta)^\beta}\widecheck\varphi*u(s)\|_2^2-2\int_s^t\langle u\cdot\nabla u(\tau),e^{-2\nu(\tau-s)(-\Delta)^\beta}\widecheck\varphi*\widecheck \varphi*u(s)\rangle\,\mathrm d\tau\\
&+2\int_s^t\langle b\cdot\nabla b(\tau),e^{-2\nu(\tau-s)(-\Delta)^\beta}\widecheck\varphi*\widecheck \varphi*u(s)\rangle\,\mathrm d\tau,
\end{split}
\end{equation}
\begin{equation}\label{u-ehigh}
\begin{split}
&E(t) \| \psi(t)\widehat u(t)\|_2^2\\
\leq& E(s)\|\psi(s)\widehat u(s)\|_2^2-2\nu\int_s^tE(\tau)\||\xi|^\beta\psi\widehat u(\tau)\|_2^2\,\mathrm d\tau+\int^t_s E'(\tau) \|\psi\widehat u(\tau)\|_2^2\, \mathrm d\tau\\
&-2\nu\int_s^tE(\tau)\langle \psi'(\tau)\widehat u(\tau),\psi\widehat u(\tau)\rangle\,\mathrm d\tau-2\int_s^tE(\tau)\langle \widehat{u\cdot\nabla u}(\tau),\psi^2\widehat u(\tau)\rangle\,\mathrm d\tau\\
&+2\int_s^tE(\tau)\langle \widehat{b\cdot\nabla b}(\tau),\psi^2\widehat u(\tau)\rangle\,\mathrm d\tau.
\end{split}
\end{equation}
\end{Lemma}

\medskip

\subsection{Auxiliary estimates}
In order to establish the decay of the high frequency parts, we need the following estimate for the Fourier transform of $u$ and $b$.

\begin{Lemma}\label{le-Fb}
Let $(u,b)$ be a mild solution to (\ref{HMHD}). Assume the initial data $(u_0,b_0)$ belongs to $\left( L^2(\R^3)\right)^2$. If $\mu>0$ and additionally $b_0\in L^1(\R^3)$, then we have
\[|\hat b(\xi, t)| \lesssim 1+\frac{1+|\xi|}{|\xi|^{2\alpha-1}}.\]
If $\nu>0$ and additionally $u_0\in L^1(\R^3)$, then we have
\[|\hat u(\xi, t)| \lesssim 1+|\xi|^{1-2\beta}.\]
\end{Lemma}
\pf
Taking Fourier transform of the second equation in (\ref{HMHD}) yields 
\[\hat b_t+|\xi|^{2\alpha}\hat b=G(\xi,t)\]
where
\[G(\xi,t)=-\widehat{u\cdot\nabla b}+\widehat{b\cdot\nabla u}-(\nabla\times((\nabla\times b)\times b))^{\wedge}.\]
Thus, we have
\[\hat b(t)=e^{-|\xi|^{2\alpha}t}\hat b(0)+\int_0^te^{-|\xi|^{2\alpha}(t-s)}G(\xi,s)\, \mathrm ds.\]
Recall the vector identity $(\nabla\times b)\times b=b\cdot\nabla b-\nabla \frac{|b|^2}2=\nabla\cdot(b\otimes b)-\nabla \frac{|b|^2}2$ and its consequence 
$\nabla\times((\nabla\times b)\times b)=\nabla\times(\nabla\cdot(b\otimes b)),$ which yields

\begin{equation}\notag
\begin{split}
|G(\xi,s)|\lesssim & \sum_{i,j}\left(|\xi| |\widehat{u^ib^j}|+|\xi|^2 |\widehat{b^ib^j}|\right)\\
\lesssim & (|\xi| \|u_0\|_2\|b_0\|_2+|\xi|^2 \|b_0\|_2^2)\\
\lesssim & |\xi|(1+|\xi|).
\end{split}
\end{equation}
It then follows that
\begin{equation}\notag
\begin{split}
|\hat b(\xi, t)|\leq & |\hat b(0)|+C|\xi|(1+|\xi|)\int_0^te^{-|\xi|^{2\alpha}(t-s)}\, \mathrm ds\\
\leq & C\|b_0\|_1+C\frac{1+|\xi|}{|\xi|^{2\alpha-1}}(1-e^{-|\xi|^{2\alpha}t}) \lesssim  1+\frac{1+|\xi|}{|\xi|^{2\alpha-1}}.
\end{split}
\end{equation}
The estimate for $\hat u$ can be established in a similar way.

\cbdu



At the end of this section, we introduce the fractional Sobolev inequality.
\begin{Lemma}\label{Sobo}\cite{MS}
Let $0\leq k<\ell \leq 1$, $1 \leq p <q< \infty$ satisfying $p(\ell-k)<n$ and $\frac{1}{q}=\frac{1}{p}-\frac{\ell-k}{n},$ then there exists a positive constant $C=C(n,p,q,k,\ell)$ such that
\[\|f\|_{W^{k,q}(\mathbb{R}^n)}\leq C\|f\|_{W^{\ell,p}(\mathbb{R}^n)}.\]
\end{Lemma}

\bigskip

\section{The inviscid resistive case  $\nu=0, \mu>0$}
\label{sec:dec}

In this section we show that for smooth solutions to the inviscid resistive Hall-MHD system (\ref{HMHD}) in $\mathbb{R}^3$, the energy $\|b(t)\|_2^2$ vanishes eventually despite the lack of velocity diffusion, provided $u(t)$ is bounded in $W^{1-\alpha,\frac3\alpha}(\mathbb{R}^3)$.
We estimate the low frequency part $\|\varphi\hat b(t)\|_2$ and high frequency part $\|\psi\hat b(t)\|_2$ separately.

\begin{Lemma}\label{le-low}
(Low frequency decay) Let $(u,b)$ be a smooth solution to (\ref{HMHD}). Assume $(u_0,b_0)\in \left(L^2(\R^3)\right)^2$. For $\varphi=e^{-|\xi|^2}$, there holds
\[\lim_{t\to \infty}\|\varphi\hat b(t)\|_2=0.\]
\end{Lemma}
\pf
The generalized energy inequality (\ref{elow}) implies
\begin{equation}\notag
\begin{split}
\|\widecheck\varphi*b(t)\|_2^2\leq &\|e^{-\mu(t-s)(-\Delta)^\alpha}\widecheck\varphi*b(t)\|_2^2\\
&+2\int_s^t\left|\langle u\cdot\nabla b(\tau),e^{-2\mu(\tau-s)(-\Delta)^\alpha}\widecheck\varphi*\widecheck \varphi*b(s)\rangle\right|\,\mathrm d\tau\\
&+2\int_s^t\left|\langle b\cdot\nabla u(\tau),e^{-2\mu(\tau-s)(-\Delta)^\alpha}\widecheck\varphi*\widecheck \varphi*b(s)\rangle\right|\,\mathrm d\tau\\
&+2\int_s^t\left|\langle\nabla\times(\nabla\cdot(b\otimes b)(\tau)), e^{-2\mu(-\Delta)^\alpha (\tau-s)}\widecheck\varphi*\widecheck\varphi* b(s)\rangle\right|\, \mathrm d\tau\\
:=&I+II+III+IV.
\end{split}
\end{equation}
One can see immediately that
 \[\limsup_{t\to \infty}I=\limsup_{t\to \infty}\|e^{-\mu |\xi|^{2\alpha}(t-s)}\varphi\hat b(s)\|_2^2=0.\]
By H\"older, Hausdorf-Young and Sobolev inequalities, the facts that $\varphi^2$ is a rapidly decreasing function and $\|u(t)\|_2$ is bounded for all the time, we infer, for $2\leq p\leq \infty$ with $\frac1p+\frac1{p'}=1$, and $\alpha\in(0,1]$
\begin{equation}\notag
\begin{split}
&\left|\langle u\cdot\nabla b(\tau), e^{-2\mu(-\Delta)^\alpha (\tau-s)}\widecheck\varphi^2* b(s)\rangle\right|\\
= & \left|\langle \xi\widehat{(u\otimes b)}(\tau), e^{-2\mu|\xi|^{2\alpha} (\tau-s)} \varphi^2\hat b(s)\rangle\right|\\ 
= & \left|\langle \widehat{(u\otimes b)}(\tau), \xi^{1-\alpha}e^{-2\mu|\xi|^{2\alpha} (\tau-s)} \varphi^2\xi^{\alpha}\hat b(s)\rangle\right|\\ 
\leq &\|\widehat{(u\otimes b)}(\tau)\|_p\|\xi^{1-\alpha} \varphi^2\|_{\frac{2p}{p-2}}\|\xi^{\alpha}e^{-2\mu|\xi|^{2\alpha} (\tau-s)}\hat b(s)\|_2\\
\lesssim &\|(u\otimes b)(\tau)\|_{p'}\|\xi^{\alpha}\hat b(\tau)\|_2\\
\lesssim & \|u(\tau)\|_2 \|b(\tau)\|_{\frac{2p}{p-2}} \|\nabla^\alpha b(\tau)\|_2\\
\lesssim &\|u(\tau)\|_2 \|\nabla^\alpha b(\tau)\|_2^2
\lesssim \|\nabla^\alpha b(\tau)\|_2^2,
\end{split}
\end{equation}
provided that $p=\frac3\alpha$ in order to apply Sobolev's inequality, which is compatible with $p\geq 2$ since $\alpha\in(0,1]$.
Applying similar strategy as above, it follows that
\begin{equation}\notag
\left|\langle b\cdot\nabla u(\tau), e^{-2\mu(-\Delta)^\alpha (\tau-s)}\widecheck\varphi^2* b(s)\rangle\right|  \lesssim \|\nabla^\alpha b(\tau)\|_2^2.
\end{equation} 
To handle the Hall term $IV$, we notice that $\||\xi|^{2-\alpha} {\varphi}^2\|_p$ is finite for any $p>1$ and $\alpha\in(0,1]$, hence
\begin{equation}\notag
\begin{split}
&\left|\langle\nabla\times(\nabla\cdot(b\otimes b))(\tau), e^{-2\mu(-\Delta)^\alpha (\tau-s)}\widecheck\varphi^2* b(s)\rangle\right|\\
=&\left|\langle\xi\times(\xi\cdot\widehat{(b\otimes b)})(\tau), e^{-2\mu|\xi|^{2\alpha} (\tau-s)}\varphi^2 \hat b(s)\rangle\right|\\
=&\left|\langle \varphi^2\xi^{-\alpha}\xi\times(\xi\cdot\widehat{(b\otimes b)})(\tau), e^{-2\mu|\xi|^{2\alpha} (\tau-s)}\xi^{\alpha} \hat b(s)\rangle\right|\\
\leq &\| |\xi|^{2-\alpha}\varphi^2\|_p \|\widehat{(b\otimes b)}(\tau)\|_{\frac{2p}{p-2}} \|\xi^{\alpha}\hat b(\tau)\|_2 \\
\lesssim & \|b(\tau)\|_2 \|b(\tau)\|_p \|\nabla^\alpha b(\tau)\|_2 \lesssim \|\nabla^\alpha b(\tau)\|_2^2
\end{split}
\end{equation} 
provided that $\alpha=\frac32-\frac3p$ and $p>2$ which is again compatible with $\alpha\in(0,1]$. 
Combining the last three inequalities yields
\[II+III+IV \lesssim \int^t_s \|\nabla^\alpha b(\tau)\|_2^2 \mathrm{d}\tau.\]
Thanks to the fact that $\int_0^\infty\|\nabla^\alpha b(t)\|_2^2\, dt$ is finite, it follows that
\[\underset{t \to \infty}{\lim} (II+III+IV) \leq \underset{s\to \infty}{\lim} \underset{t\to \infty}{\lim} \int^t_s \|\nabla^\alpha b(\tau)\|^2_2 \mathrm{d}\tau =0.\]
Therefore, we conclude \[\lim_{t\to \infty}\|\varphi\widehat b(t)\|_2=\lim_{t\to \infty}\|\widecheck\varphi*b(t)\|_2=0.\]

\cbdu

\begin{Lemma} \label{le-high}
(High frequency decay)Let $(u,b)$ be a smooth solution to (\ref{HMHD}). Assume $u_0\in L^2(\R^3)$ and $b_0\in L^1(\R^3)\cap L^2(\R^3)$, and one of the two conditions holds: 
\begin{equation}\notag
\begin{split}
(i)\ \ &  u\in L^\infty(0,\infty;W^{1-\alpha, \frac3\alpha}(\R^3)) \text{ with }  \alpha\in[1/2, 1]; \\
(ii)\ \ & b \in L^\infty(0,\infty;W^{1-\alpha,\infty}(\R^3)) \text{ and } u\in L^\infty(0,\infty;W^{1-\alpha, \frac3\alpha}(\R^3)) \text{ with } \alpha\in (0,1].
\end{split}
\end{equation}
Then 
\[\lim_{t\to 0}\|\psi\hat b(t)\|_2=0.\]
\end{Lemma}
\pf
We start with estimating the last three integrals on the right hand side of the generalized energy inequality (\ref{ehigh}), recalled here,
\begin{equation}\notag
\begin{split}
E(t) \| \psi\widehat b(t)\|_2^2\leq& E(s)\|\psi\widehat b(s)\|_2^2-2\mu\int_s^tE(\tau)\||\xi|^\alpha\psi\widehat b(\tau)\|_2^2\,\mathrm d\tau
+\int^t_s E'(\tau) \|\psi\widehat b(\tau)\|_2^2\, \mathrm d\tau\\
&+2\int_s^tE(\tau)\langle \widehat{b\cdot\nabla u}(\tau),\psi^2\widehat b(\tau)\rangle\,\mathrm d\tau
-2\int_s^tE(\tau)\langle \widehat{u\cdot\nabla b}(\tau),\psi^2\widehat b(\tau)\rangle\,\mathrm d\tau\\
&-2\int^t_s E(\tau)\langle{\nabla\times(\nabla\cdot(b\times b)(\tau))}^\wedge , \psi^2\widehat b(\tau)\rangle\, \mathrm d\tau\\
:=&J_0+J_1+J_2+J_3+J_4+J_5.
\end{split}
\end{equation}
In order to estimate $J_3$ where no cancelation presents, we need the additional assumptions on $u$ and $b$. First, we have by using H\"older's inequality and Plancherel's theorem 
\begin{equation}\notag
\begin{split}
&\int^t_s E(\tau)\left|\langle\widehat{b\cdot \nabla u}(\tau), \psi^2 \hat{b}(\tau)\rangle\right|\mathrm{d}\tau \\
= & \int^t_s E(\tau)\left|\langle\xi^{-\alpha}\xi\cdot \widehat{b\otimes u}(\tau),\psi^2 \xi^{\alpha}\hat{b}(\tau)\rangle\right|\mathrm{d}\tau\\
\leq & \int^t_s E(\tau)\||\xi|^{1-\alpha} \widehat{b\otimes u}(\tau)\|_2\|\psi^2 \xi^{\alpha}\hat b(\tau)\|_2\mathrm{d}\tau\\
\lesssim & \int^t_s E(\tau)\left(\|u\nabla^{1-\alpha}b(\tau)\|_2+\|b\nabla^{1-\alpha}u(\tau)\|_2\right)\|\nabla^{\alpha} b(\tau)\|_2\mathrm{d}\tau.
\end{split}
\end{equation}
If $u\in L^\infty(0,\infty; W^{1-\alpha, \frac3\alpha}(\mathbb{R}^3))$, it follows from H\"older's inequality and the fractional Sobolev inequality in Lemma \ref{Sobo} that
\begin{equation}\notag
\|b\nabla^{1-\alpha}u(\tau)\|_2\leq \|b\|_{\frac6{3-2\alpha}}\|\nabla^{1-\alpha}u\|_{\frac3\alpha}
\lesssim \|\nabla^\alpha b\|_2\|\nabla^{1-\alpha}u\|_{\frac3\alpha}.
\end{equation}
If $\alpha\geq\frac12$, then it follows from H\"older's inequality and the Sobolev inequality that
\begin{equation}\notag
\|u\nabla^{1-\alpha}b(\tau)\|_2\leq \|u\|_{\frac3{2\alpha-1}}\|\nabla^{1-\alpha}b\|_{\frac6{5-4\alpha}}
\lesssim \|\nabla^{1-\alpha}u\|_{\frac3\alpha}\|\nabla^{\alpha}b\|_2.
\end{equation}
While if $b\in L^\infty(0,\infty; W^{1-\alpha,\infty}(\mathbb{R}^3))$, $\|u\nabla^{1-\alpha}b(\tau)\|_2$ can be also estimated as 
\[\|u\nabla^{1-\alpha}b(\tau)\|_2\leq \|u(\tau)\|_2\|\nabla^{1-\alpha}b(\tau)\|_\infty.\]
As a conclusion, if $u\in L^\infty(0,\infty; W^{1-\alpha, \frac3\alpha}(\mathbb{R}^3))$ with $\frac12\leq \alpha\leq 1$, we have 
\begin{equation}\label{j3-2}
\begin{split}
&\int^t_s E(\tau)\langle\widehat{b\cdot \nabla u}(\tau), \psi^2 \hat{b}(\tau)\rangle\mathrm{d}\tau \\
\lesssim & \int^t_s E(\tau) \|\nabla^{1-\alpha}u\|_{\frac3\alpha}\|\nabla^\alpha b(\tau)\|_2^2\mathrm{d}\tau
\lesssim  \int^t_s E(\tau)\|\nabla^\alpha b(\tau)\|^2_2\mathrm{d}\tau;
\end{split}
\end{equation}
If $b \in L^\infty(0,\infty; W^{1-\alpha,\infty}(\mathbb{R}^3))$ and $u\in L^\infty(0,\infty;W^{1-\alpha, \frac3\alpha}(\mathbb{R}^3))$ with $\alpha\in (0,1]$, we have
\begin{equation}\label{j3-1}
\begin{split}
&\int^t_s E(\tau)\langle\widehat{b\cdot \nabla u}(\tau), \psi^2 \hat{b}(\tau)\rangle\mathrm{d}\tau \\
\lesssim& \int^t_s E(\tau)(\|\nabla^\alpha b\|_2\|\nabla^{1-\alpha}u\|_{\frac3\alpha}+\|u(\tau)\|_2\|\nabla^{1-\alpha}b(\tau)\|_\infty)\|{\nabla^\alpha b(\tau)}\|_2\mathrm{d}\tau\\
\lesssim & \Big(\int^t_s E(\tau)^2 \mathrm{d}\tau\Big)^{\frac{1}{2}} \Big(\int^t_s \|\nabla b(\tau)\|^2_2 \mathrm{d}\tau\Big)^{\frac{1}{2}}+\int^t_s E(\tau)\|\nabla^\alpha b(\tau)\|^2_2\mathrm{d}\tau.
\end{split}
\end{equation}
To deal with $J_4$, observing the cancelation $\langle u\cdot \nabla b, b\rangle=0$ 
we have, 
\begin{equation}\notag
\begin{split}
\int^t_s E(\tau)\left|\langle \widehat{u\cdot \nabla b}(\tau), \psi^2 \hat{b}(\tau)\rangle\right| \mathrm{d}\tau= &\int^t_s E(\tau)\left|\langle\widehat{u\cdot \nabla b}(\tau), (\psi^2-1)\hat{b}(\tau)\rangle\right|\mathrm{d}\tau\\
= &\int^t_s E(\tau)\left|\langle \widehat{u\otimes  b}(\tau), (\psi^2-1)\widehat{\nabla b}(\tau)\rangle\right|\mathrm{d}\tau\\
\leq &\int^t_s E(\tau)\|\widehat{u\otimes  b}(\tau)\|_{\frac3\alpha} \|\widehat{\nabla^{\alpha} b}(\tau)\|_2\||\xi|^{1-\alpha}(\psi^2-1)\|_{\frac6{3-2\alpha}}\mathrm{d}\tau
\end{split}
\end{equation}
Notice that $\psi^2-1=-2e^{-|\xi|^2}+e^{-2|\xi|^2}$ and $\||\xi|^{1-\alpha}(\psi^2-1)\|_p$ is finite for any $p>1$ and $\alpha\leq 1$, we continue the estimate as
\begin{equation}\notag
\begin{split}
\int^t_s E(\tau)\left|\langle \widehat{u\cdot \nabla b}(\tau), \psi^2 \hat{b}(\tau)\rangle\right| \mathrm{d}\tau
\lesssim &\int^t_s E(\tau)\|u(\tau)\|_2\|b(\tau)\|_{\frac6{3-2\alpha}}\|\nabla^{\alpha}b(\tau)\|_2 \mathrm{d}\tau\\
\lesssim & \int^t_s E(\tau) \|\nabla^\alpha b(\tau)\|_2^2 \mathrm{d}\tau
\end{split}
\end{equation}
where we apply Sobolev's inequality to obtain the last step.
Again, utilizing the cancelation $\langle \nabla \times ((\nabla \times b)\times b),b\rangle=0$ and the fact $\||\xi|^{2-\alpha}(\psi^2-1)\|_p$ is finite for $p>1$, $J_5$ can be estimated as
\begin{equation}\notag
\begin{split}
&\int^t_s E(\tau)\left|\langle[\nabla \times(\nabla \cdot (b \otimes b))]^{\wedge}, \psi^2 \hat{b}\rangle (\tau)\right|\mathrm{d}\tau \\
= &\int^t_s E(\tau)\left|\langle[\nabla \times(\nabla \cdot (b \otimes b))]^{\wedge}, (\psi^2-1) \hat{b}\rangle (\tau)\right|\mathrm{d}\tau\\
= &\int^t_s E(\tau)\left|\langle (\psi^2-1) \xi^{-\alpha}[\nabla \times(\nabla \cdot (b \otimes b))]^{\wedge},  \xi^{\alpha}\hat{b}\rangle (\tau)\right|\mathrm{d}\tau\\
\lesssim & \int^t_s E(\tau)\||\xi|^{2-\alpha}(\psi^2-1)\|_{\frac6{3-2\alpha}}\|\widehat{(b\otimes b)}(\tau)\|_{\frac3\alpha}\|\widehat{\nabla^\alpha b}(\tau)\|_2 \mathrm{d}\tau\\
\lesssim & \int^t_s E(\tau)\|b(\tau)\|_2\|b(\tau)\|_{\frac6{3-2\alpha}}\|\nabla^\alpha b(\tau)\|_2 \mathrm{d}\tau\\
\lesssim & \int^t_s E(\tau) \|\nabla^\alpha b(\tau)\|_2^2\mathrm{d}\tau.
\end{split}
\end{equation}
If $b \in L^\infty(0,\infty; W^{1-\alpha,\infty}(\mathbb{R}^3))$ and $u\in L^\infty(0,\infty;W^{1-\alpha, \frac3\alpha}(\mathbb{R}^3))$ with $\alpha\in (0,1]$, combining (\ref{ehigh}), (\ref{j3-1}) and the last three inequalities yields
\begin{equation}\label{ehigh-1}
\begin{split}
\| \psi\widehat b(t)\|_2^2\leq& \frac{E(s)}{E(t)}\|\psi\widehat b(s)\|_2^2-2\mu\int_s^t\frac{E(\tau)}{E(t)}\||\xi|^{\alpha}\psi\widehat b(\tau)\|_2^2\,\mathrm d\tau
+\int^t_s \frac{E'(\tau)}{E(t)} \|\psi\widehat b(\tau)\|_2^2\, \mathrm d\tau\\
&+ \int^t_s \frac{E(\tau)}{E(t)} \|\nabla^\alpha b(\tau)\|_2^2\mathrm{d}\tau+\frac1{E(t)}\Big(\int^t_s E(\tau)^2 \mathrm{d}\tau\Big)^{\frac{1}{2}} \Big(\int^t_s \|\nabla^\alpha b(\tau)\|^2_2 \mathrm{d}\tau\Big)^{\frac{1}{2}};
\end{split}
\end{equation}
If $u\in L^\infty(0,\infty; W^{1-\alpha, \frac3\alpha}(\mathbb{R}^3))$ with $\frac12\leq \alpha\leq 1$, combining (\ref{ehigh}), (\ref{j3-2}) and the these estimates for $J_3, J_4, J_5$ yields
\begin{equation}\notag
\begin{split}
\| \psi\widehat b(t)\|_2^2\leq& \frac{E(s)}{E(t)}\|\psi\widehat b(s)\|_2^2-2\mu\int_s^t\frac{E(\tau)}{E(t)}\||\xi|^{\alpha}\psi\widehat b(\tau)\|_2^2\,\mathrm d\tau\\
&+\int^t_s \frac{E'(\tau)}{E(t)} \|\psi\widehat b(\tau)\|_2^2\, \mathrm d\tau
+ \int^t_s \frac{E(\tau)}{E(t)} \|\nabla^\alpha b(\tau)\|_2^2\mathrm{d}\tau.
\end{split}
\end{equation}
We only need to take care of (\ref{ehigh-1}) for the obvious reason. We will apply Fourier splitting approach to handle the second and third terms on the right hand of (\ref{ehigh-1}). Denote the ball $B(t)=\{\xi \in \mathbb{R}^3: |\xi| \leq G(t)\}$, where the radius $G(t)$ will be determined later. Then we infer
\begin{equation}\notag
\begin{split}
&-2\mu\int_s^t\frac{E(\tau)}{E(t)}\||\xi|^{\alpha}\psi\widehat b(\tau)\|_2^2\,\mathrm d\tau
+\int^t_s \frac{E'(\tau)}{E(t)} \|\psi\widehat b(\tau)\|_2^2\, \mathrm d\tau\\
\leq&-2\mu\int_s^t\frac{E(\tau)}{E(t)}\int_{B^c(t)}||\xi|^{\alpha}\psi\widehat b(\tau)|^2\, \mathrm d\xi \,\mathrm d\tau+\int^t_s \frac{E'(\tau)}{E(t)} \int_{B^c(t)}|\psi\widehat b(\tau)|^2\, \mathrm d\xi\, \mathrm d\tau\\
&+\int^t_s \frac{E'(\tau)}{E(t)} \int_{B(t)}|\psi\widehat b(\tau)|^2\, \mathrm d\xi\, \mathrm d\tau\\
\leq & \int^t_s \frac{E'(\tau)-2\mu E(\tau)G^{2\alpha}(\tau)}{E(t)} \int_{B^c(t)}|\psi\widehat b(\tau)|^2\, \mathrm d\xi\, \mathrm d\tau
+\int^t_s \frac{E'(\tau)}{E(t)} \int_{B(t)}|\psi\widehat b(\tau)|^2\, \mathrm d\xi\, \mathrm d\tau\\
\end{split}
\end{equation}
Take \[E(t)=e^{\varepsilon t},\text{ and }G(t)=(\frac{\varepsilon}{2\mu})^{\frac1{2\alpha}},\] 
which indicates that $E'(t)-2\mu E(t)G^{2\alpha}(t)=0.$ Therefore, it follows from (\ref{ehigh-1})
\begin{equation}\notag
\begin{split}
\| \psi\widehat b(t)\|_2^2\leq& \frac{E(s)}{E(t)}\|\psi\widehat b(s)\|_2^2+\int^t_s \frac{E'(\tau)}{E(t)} \int_{B(t)}|\psi\widehat b(\tau)|^2\, \mathrm d\xi\, \mathrm d\tau\\
&+ \int^t_s \frac{E(\tau)}{E(t)} \|\nabla^\alpha b(\tau)\|_2^2\mathrm{d}\tau+\frac1{E(t)}\Big(\int^t_s E(\tau)^2 \mathrm{d}\tau\Big)^{\frac{1}{2}} \Big(\int^t_s \|\nabla^\alpha b(\tau)\|^2_2 \mathrm{d}\tau\Big)^{\frac{1}{2}};
\end{split}
\end{equation}
By Lemma \ref{le-Fb} we can bound the second term,  since $\psi^2 \leq 1,$
\begin{equation}\notag
\begin{split}
\int_{B(t)}\left| \psi \hat{b}(\tau)\right|^2\mathrm{d}\xi \lesssim & \int_{B(t)} \Big(1+\frac{1}{|\xi|^{2\alpha-1}}\Big)^2\mathrm{d}\xi\lesssim  \int_{B(t)} \Big(1+\frac{1}{|\xi|^{4\alpha-2}}\Big)\mathrm{d}\xi\\
\lesssim &\int_0^{G(t)}(1+r^{2-4\alpha})r^2 \mathrm{d}r
\lesssim \varepsilon^{\frac{3}{2\alpha}}+\varepsilon^{\frac{5-4\alpha}{2\alpha}}.
\end{split}
\end{equation}
Thus, the estimate is continued as 
\begin{equation}\notag
\begin{split}
\|\psi \hat b(t)\|_2^2 \leq & \frac{E(s)}{E(t)}\|\psi \hat{b}(s)\|_2^2 +\frac{C}{e^{\varepsilon t}}\Big(\frac{e^{2\varepsilon t}}{2 \varepsilon}\Big)^{\frac{1}{2}} \Big(\int^t_s \|\nabla^\alpha b(\tau)\|^2_2 \mathrm{d}\tau\Big)^{\frac{1}{2}}\\
& +\frac{C}{E(t)}\int^t_s E(\tau)\|\nabla^\alpha b(\tau)\|^2_2\mathrm{d}\tau+C(\varepsilon^{\frac{3}{2\alpha}}+\varepsilon^{\frac{5-4\alpha}{2\alpha}})
\end{split}
\end{equation} 
with various constants $C$ which are independent of $t,s$ and $\varepsilon$.
Now, we first pass the limit $t\to \infty$, 
\begin{equation}\notag
\begin{split}
\underset{{t \to \infty}}{\lim}{\|\psi \hat b(t)\|_2^2} \leq & \underset{{t \to \infty}}{\lim}\frac{E(s)}{E(t)}\|\psi \hat{b}(s)\|_2^2 +\underset{{t \to \infty}}{\lim}\frac{C}{\sqrt{\varepsilon}}\Big(\int^t_s \|\nabla^\alpha b(\tau)\|^2_2 \mathrm{d}\tau\Big)^{\frac{1}{2}}\\
&+\underset{{t \to \infty}}{\lim}\frac{C}{E(t)}\int^t_s E(\tau)\|\nabla^\alpha b(\tau)\|^2_2\mathrm{d}\tau
+C(\varepsilon^{\frac{3}{2\alpha}}+\varepsilon^{\frac{5-4\alpha}{2\alpha}})\\
\leq & \underset{{t \to \infty}}{\lim}e^{\varepsilon (s-t)}\|b_0\|_2^2+\frac{C}{\sqrt{\varepsilon}}\Big(\int^\infty_s \|\nabla^\alpha b(\tau)\|^2_2 \mathrm{d}\tau\Big)^{\frac{1}{2}}\\
&+ C\int^\infty_s \|\nabla^\alpha b(\tau)\|^2_2\mathrm{d}\tau +C(\varepsilon^{\frac{3}{2\alpha}}+\varepsilon^{\frac{5-4\alpha}{2\alpha}})\\
\leq &C(\varepsilon^{\frac{3}{2\alpha}}+\varepsilon^{\frac{5-4\alpha}{2\alpha}})+\frac{C}{\sqrt{\varepsilon}}\Big(\int^\infty_s \|\nabla^\alpha b(\tau)\|^2_2 \mathrm{d}\tau\Big)^{\frac{1}{2}}+ C\int^\infty_s \|\nabla^\alpha b(\tau)\|^2_2\mathrm{d}\tau;
\end{split}
\end{equation} 
and then pass the limit $s\to \infty$,
\begin{equation}\notag
\begin{split}
&\underset{{t \to \infty}}{\lim}{\|\psi \hat b(t)\|_2^2} \\
\leq &\underset{{s \to \infty}}{\lim}\Bigg(C(\varepsilon^{\frac{3}{2\alpha}}+\varepsilon^{\frac{5-4\alpha}{2\alpha}})+\frac{C}{\sqrt{\varepsilon}}\Big(\int^\infty_s \|\nabla^\alpha b(\tau)\|^2_2 \mathrm{d}\tau\Big)^{\frac{1}{2}}+C\int^\infty_s \|\nabla^\alpha b(\tau)\|^2_2\mathrm{d}\tau\Bigg)\\
\leq &C(\varepsilon^{\frac{3}{2\alpha}}+\varepsilon^{\frac{5-4\alpha}{2\alpha}}).
\end{split}
\end{equation}
Since $\varepsilon>0$ can be chosen arbitrarily small, it implies that $\underset{{t \to \infty}}{\lim}\|\psi \hat b(t)\|_2=0.$

\cbdu

{\bf {Proof of Theorem \ref{thm1}:}} Combining Lemma \ref{le-low} and Lemma \ref{le-high} yields 
\[\lim _{t\to\infty}\|b(t)\|_2^2=0.\]
This convergence along with the basic energy law implies $\underset{{t \to \infty}}{\lim}\| u(t)\|_2^2=C$ for a constant $C$.

\section{The viscous non-resistive case $\nu>0, \mu=0$}
\label{sec:dec1}

In this section we show that $\|u(t)\|_2$ converges to zero in the viscous non-resistive case $\nu>0,\mu=0$ and $\beta>0$, provided $b$ is bounded in $W^{1-\alpha,\infty}(\R^3).$
We estimate the low frequency part $\|\varphi(t)\hat u(t)\|_2$ and high frequency part $\|(1-\varphi(t))\hat u(t)\|_2$ separately, by taking $\varphi(t)=e^{-|\xi|^{2\beta}t}$.

\begin{Lemma}(Low frequency decay) Let $(u,b)$ be a smooth solution to (\ref{HMHD}) with $\nu>0$ and $\mu=0$. Assume $(u_0,b_0)\in \left(L^2(\R^3)\right)^2$. For $\varphi(t)=e^{-|\xi|^{2\beta}t}$, there holds
\[\lim_{t\to \infty}\|\varphi(t)\hat u(t)\|_2=0.\]
\end{Lemma}
\pf
The generalized energy inequality (\ref{u-elow}) implies
\begin{equation}\notag
\begin{split}
\|\widecheck\varphi*u(t)\|_2^2\leq &\|e^{-\nu(t-s)(-\Delta)^\beta}\widecheck\varphi*u(t)\|_2^2\\
&+2\int_s^t\left|\langle u\cdot\nabla u(\tau),e^{-2\nu(\tau-s)(-\Delta)^\beta}\widecheck\varphi*\widecheck \varphi*u(s)\rangle\right|\,\mathrm d\tau\\
&+2\int_s^t\left|\langle b\cdot\nabla b(\tau),e^{-2\nu(\tau-s)(-\Delta)^\beta}\widecheck\varphi*\widecheck \varphi*u(s)\rangle\right|\,\mathrm d\tau\\
:=&I+II+III.
\end{split}
\end{equation}
One can see immediately that
 \[\limsup_{t\to \infty}I=\limsup_{t\to \infty}\|e^{-\nu |\xi|^{2\beta}(t-s)}\varphi(s)\hat u(s)\|_2^2=0.\]
By H\"older, Hausdorf-Young and Sobolev inequalities, the facts that $\varphi^2(t)$ is a rapidly decreasing function of $|\xi|$ and $\|u(t)\|_2$ is bounded for all the time, we infer, for $\beta\in(0,1]$
\begin{equation}\notag
\begin{split}
&\left|\langle u\cdot\nabla u(\tau), e^{-2\nu(-\Delta)^\beta (\tau-s)}\widecheck\varphi^2* u(s)\rangle\right| \\
= & \left|\langle \xi\widehat{(u\otimes u)}(\tau), e^{-2\nu|\xi|^{2\beta} (\tau-s)} \varphi^2\hat u(s)\rangle\right|\\ 
= & \left|\langle \widehat{(u\otimes u)}(\tau), \xi^{1-\beta}e^{-2\nu|\xi|^{2\beta} (\tau-s)} \varphi^2\xi^{\alpha}\hat u(s)\rangle\right|\\ 
\leq &\|\widehat{(u\otimes u)}(\tau)\|_{\frac3\beta}\|\xi^{1-\beta} \varphi^2\|_{\frac{6}{3-2\beta}}\|\xi^{\beta}e^{-2\nu|\xi|^{2\beta} (\tau-s)}\hat u(s)\|_2\\
\lesssim &\|(u\otimes u)(\tau)\|_{\frac3{3-\beta}}\|\xi^{\beta}\hat u(\tau)\|_2\\
\lesssim & \|u(\tau)\|_2 \|u(\tau)\|_{\frac{6}{3-2\beta}} \|\nabla^\beta u(\tau)\|_2\\
\lesssim &\|u(\tau)\|_2 \|\nabla^\beta u(\tau)\|_2^2
\lesssim \|\nabla^\beta u(\tau)\|_2^2.
\end{split}
\end{equation}
Applying similar strategy as above, it follows that
\begin{equation}\notag
\begin{split}
&\left|\langle b\cdot\nabla b(\tau), e^{-2\nu(-\Delta)^\beta (\tau-s)}\widecheck\varphi^2* u(s)\rangle\right| \\
= & \left|\langle \xi^{1-\beta}\varphi^2\widehat{(b\otimes b)}(\tau), e^{-2\nu|\xi|^{2\beta} (\tau-s)}\xi^{\beta}\hat u(s)\rangle\right|\\ 
\leq &\|\widehat{(b\otimes b)}(\tau)\|_{\infty}\|\xi \varphi^2(\tau)\|_2\|e^{-2\nu|\xi|^{2\beta} (\tau-s)}\hat u(s)\|_2\\
\lesssim &\|(b\otimes b)(\tau)\|_{1}\|\xi \varphi^2(\tau)\|_2\|\hat u(\tau)\|_2\\
\lesssim & \|b(\tau)\|_2^2\|\xi\varphi^2(\tau)\|_2 \|u(\tau)\|_2
\lesssim  \|\xi\varphi^2(\tau)\|_2.
\end{split}
\end{equation}
Combining the last two inequalities yields
\[II+III \lesssim \int^t_s \|\nabla^\beta u(\tau)\|_2^2 \mathrm{d}\tau+\int^t_s \|\xi\varphi^2(\tau)\|_2 \mathrm{d}\tau.\]
Straightforward computation shows that
\begin{equation}\notag
\begin{split}
\|\xi\varphi^2(\tau)\|_2^2=&\int_{\R^3}|\xi|^2e^{-4|\xi|^{2\beta}\tau}\, \mathrm d\xi\lesssim 
\int_0^\infty r^4e^{-4r^{2\beta}\tau}\, \mathrm dr\\
=&\tau^{-\frac5{2\beta}} \int_0^\infty w^4e^{-4w^{2\beta}}\, \mathrm dw\lesssim \tau^{-\frac5{2\beta}}.
\end{split}
\end{equation}
Together with the fact that $\int_0^\infty\|\nabla^\beta u(t)\|_2^2\, dt$ is finite due to the basic energy law, it follows 
\[\underset{t \to \infty}{\lim} (II+III) \leq \underset{s\to \infty}{\lim} \underset{t\to \infty}{\lim} \int^t_s \left(\|\nabla^\beta u(\tau)\|^2_2+ \tau^{-\frac5{4\beta}}\right)\, \mathrm{d}\tau =0, \ \ \mbox { for } 0<\beta\leq 1.\]
Therefore, we conclude \[\lim_{t\to \infty}\|\varphi(t)\widehat u(t)\|_2=\lim_{t\to \infty}\|\widecheck\varphi(t)*u(t)\|_2=0.\]

\cbdu

\begin{Lemma} (High frequency decay)Let $(u,b)$ be a smooth solution to (\ref{HMHD}). Assume that  $u_0 \in L^1(\R^3)\cap L^2(\R^3)$, $b_0 \in L^2(\R^3)$ and $b \in L^{\infty}(0,\infty; W^{1-\alpha,\infty}(\R^3))$, then for $\psi(t)=1-e^{-|\xi|^{2\beta}t}$
\[\lim_{t\to 0}\|\psi(t)\hat u(t)\|_2=0.\]
\end{Lemma}
\pf
We start with estimating the last three integrals on the right hand side of the generalized energy inequality (\ref{u-ehigh}), recalled here,
\begin{equation}\notag
\begin{split}
&E(t) \| \psi(t)\widehat u(t)\|_2^2\\
\leq& E(s)\|\psi(s)\widehat u(s)\|_2^2-2\nu\int_s^tE(\tau)\||\xi|^\beta\psi(\tau)\widehat u(\tau)\|_2^2\,\mathrm d\tau
+\int^t_s E'(\tau) \|\psi(\tau)\widehat u(\tau)\|_2^2\, \mathrm d\tau\\
&-2\nu\int_s^tE(\tau)\langle \psi'(\tau)\widehat u(\tau),\psi(\tau)\widehat u(\tau)\rangle\,\mathrm d\tau+2\int_s^tE(\tau)\langle \widehat{b\cdot\nabla b}(\tau),\psi^2(\tau)\widehat u(\tau)\rangle\,\mathrm d\tau\\
&-2\int_s^tE(\tau)\langle \widehat{u\cdot\nabla u}(\tau),\psi^2(\tau)\widehat u(\tau)\rangle\,\mathrm d\tau\\
:=&J_0+J_1+J_2+J_3+J_4+J_5.
\end{split}
\end{equation}
Recalling that $\psi(t)=1-e^{-|\xi|^{2\beta}t}$ and $\psi'=|\xi|^{2\beta}\varphi$, $J_3$ can be estimated as
\begin{equation}\label{uj3}
\begin{split}
&\int_s^tE(\tau)\left|\langle \psi'(\tau)\widehat u(\tau),\psi(\tau)\widehat u(\tau)\rangle\right|\,\mathrm d\tau\\
=&\int_s^tE(\tau)\left|\langle |\xi|^{2\beta}\varphi(\tau)\widehat u(\tau),\psi(\tau)\widehat u(\tau)\rangle\right|\,\mathrm d\tau\\
\leq&\int_s^tE(\tau)\|\nabla^\beta u(\tau)\|_2^2\,\mathrm d\tau.
\end{split}
\end{equation}
In order to estimate $J_4$ where no cancelation presents, we need the additional assumptions on $u$ and $b$. First, we have by using H\"older's inequality and Plancherel's theorem 
\begin{equation}\label{uj4}
\begin{split}
\int^t_s E(\tau)\left|\langle\widehat{b\cdot \nabla b}(\tau), \psi^2 \hat{u}(\tau)\rangle\right|\mathrm{d}\tau 
= & \int^t_s E(\tau)\left|\langle\xi^{-\beta}\xi\cdot \widehat{b\otimes b}(\tau),\psi^2 \xi^{\beta}\hat{u}(\tau)\rangle\right|\mathrm{d}\tau\\
\leq & \int^t_s E(\tau)\||\xi|^{1-\beta} \widehat{b\otimes b}(\tau)\|_2\|\psi^2 \xi^{\beta}\hat u(\tau)\|_2\mathrm{d}\tau\\
\lesssim & \int^t_s E(\tau)\|b\nabla^{1-\beta}b(\tau)\|_2\|\nabla^{\beta} u(\tau)\|_2\mathrm{d}\tau\\
\lesssim & \int^t_s E(\tau)\|b\|_2\|\nabla^{1-\beta}b(\tau)\|_\infty\|\nabla^{\beta} u(\tau)\|_2\mathrm{d}\tau\\
\lesssim & \int^t_s E(\tau)\|\nabla^{\beta} u(\tau)\|_2\mathrm{d}\tau
\end{split}
\end{equation}
provided $b\in L^\infty (0,\infty; W^{1-\beta, \infty}(\mathbb{R}^3))$.
To deal with $J_5$, observing the cancelation $\langle u\cdot \nabla u, u\rangle=0$ 
we have, 
\begin{equation}\notag
\begin{split}
\int^t_s E(\tau)\left|\langle \widehat{u\cdot \nabla u}(\tau), \psi^2 \hat{u}(\tau)\rangle\right| \mathrm{d}\tau= &\int^t_s E(\tau)\left|\langle\widehat{u\cdot \nabla u}(\tau), (\psi^2-1)\hat{u}(\tau)\rangle\right|\mathrm{d}\tau\\
= &\int^t_s E(\tau)\left|\langle \widehat{u\otimes  u}(\tau), (\psi^2-1)\widehat{\nabla u}(\tau)\rangle\right|\mathrm{d}\tau\\
\leq &\int^t_s E(\tau)\|\widehat{u\otimes  u}(\tau)\|_{\frac3\beta} \|\widehat{\nabla^{\beta} u}(\tau)\|_2\||\xi|^{1-\beta}(\psi^2-1)\|_{\frac6{3-2\beta}}\mathrm{d}\tau
\end{split}
\end{equation}
Notice that $\psi^2-1=-2e^{-|\xi|^{2\beta}t}+e^{-2|\xi|^{2\beta}t}$ and $\||\xi|^{1-\beta}(\psi^2-1)\|_p$ is finite for any $p>1$ and $\beta\leq 1$, we continue the estimate as
\begin{equation}\label{uj5}
\begin{split}
\int^t_s E(\tau)\left|\langle \widehat{u\cdot \nabla u}(\tau), \psi^2 \hat{u}(\tau)\rangle\right| \mathrm{d}\tau
\lesssim &\int^t_s E(\tau)\|u(\tau)\|_2\|u(\tau)\|_{\frac6{3-2\beta}}\|\nabla^{\beta}u(\tau)\|_2 \mathrm{d}\tau\\
\lesssim & \int^t_s E(\tau) \|\nabla^\beta u(\tau)\|_2^2 \mathrm{d}\tau
\end{split}
\end{equation}
where we apply Sobolev's inequality to obtain the last step.

Combining (\ref{u-ehigh}) and (\ref{uj3})-(\ref{uj5}) yields
\begin{equation}\notag
\begin{split}
\| \psi\widehat u(t)\|_2^2\leq& \frac{E(s)}{E(t)}\|\psi\widehat u(s)\|_2^2-2\nu\int_s^t\frac{E(\tau)}{E(t)}\||\xi|^{\beta}\psi\widehat u(\tau)\|_2^2\,\mathrm d\tau
+\int^t_s \frac{E'(\tau)}{E(t)} \|\psi\widehat u(\tau)\|_2^2\, \mathrm d\tau\\
&+ \int^t_s \frac{E(\tau)}{E(t)} \|\nabla^\beta u(\tau)\|_2^2\mathrm{d}\tau+\frac1{E(t)}\Big(\int^t_s E(\tau)^2 \mathrm{d}\tau\Big)^{\frac{1}{2}} \Big(\int^t_s \|\nabla^\beta u(\tau)\|^2_2 \mathrm{d}\tau\Big)^{\frac{1}{2}}
\end{split}
\end{equation}
which has the same form of (\ref{ehigh-1}).
Therefore, we can apply the same Fourier splitting strategy to obtain that $\underset{{t \to \infty}}{\lim}\|\psi \hat u(t)\|_2=0.$
\cbdu

Thus the statement of Theorem \ref{thm2} follows from the two lemmas above and the basic energy equality.

\bigskip

\medskip



\bigskip



\begin{thebibliography}{XX}



\bibitem{ADFL}
M. Acheritogaray, P. Degond, A. Frouvelle and J-G. Liu.
\newblock {\em Kinetic formulation and global existence for the Hall-Magnetohydrodynamic system}.
\newblock Kinetic and Related Models, 4: 901--918, 2011.

\bibitem{AS}
R. Agapito and M. E. Schonbek.
\newblock{\em Non-uniform decay of MHD equations with and without magnetic diffusion}.
\newblock Comm. Partial Differential Equations 32, 1791-1812, 2007.



\bibitem{CDL}
D. Chae, P. Degond and J-G. Liu.
\newblock {\em Well-posedness for Hall--magnetohydrodynamics}.
\newblock arXiv:1212.3919, 2012.

\bibitem{CL}
D. Chae and J. Lee.
\newblock {\em On the blow-up criterion and small data global existence for the Hall-magneto-hydrodynamics}.
\newblock J. Differential Equations, 256: 3835--3858, 2014.

\bibitem{CS}
D. Chae,  and M. Schonbek.
\newblock {\em On the temporal decay for the Hall-magnetohydrodynamic equations}.
\newblock J. Differ. Equ., 255 (11): 3971--3982, 2013.

\bibitem{CWW}
D. Chae,  R. Wan and J. Wu.
\newblock {\em Local well-posedness for the Hall--MHD equations with fractional magnetic diffusion}.
\newblock arXiv:1404.0486v2, 2014.








\bibitem{Dai}
M. Dai.
\newblock {\em Regularity criterion for the 3D Hall-magneto-hydrodynamics}.
\newblock  Journal of Differential Equations, Vol. 261: 573--591, 2016.

\bibitem{DFRSS}
M. Dai, E. Feireisl, E. Rocca, G. Schimperna, and M. E. Schonbek.
\newblock {\em On asymptotic isotropy for a hydrodynamic model of liquid crystals}.
\newblock  Asymptotic Analysis, Vol. 97 (3-4): 189--210, 2016.

\bibitem{DQS}
M. Dai, J. Qing and M. E. Schonbek.
\newblock {\em Asymptotic behavior of solutions to liquid crystal systems in $\mathbb{R}^3$}.
\newblock  Communications in Partial Differential Equations, Vol. 37(12): 2138--2164, 2012.

\bibitem{DSch}
M. Dai and M. Schonbek.
\newblock {\em Asymptotic behavior of solutions to the liquid crystal systems in $H^m(\mathbb{R}^3)$}.
\newblock  SIAM Journal on Mathematical Analysis, Vol. 46(5): 3131--3150, 2014.

\bibitem{DS}
E. Dumas and F. Sueur.
\newblock {\em On the weak solutions to the Maxwell-Landau-Lifshitz equations and to the Hall-magnetohydrodynamic equations}.
\newblock Comm. Math. Phys., 330: 1179--1225, 2014.



\bibitem{For}
T. G. Forbes.
\newblock {\em Magnetic reconnection in solar flares}.
\newblock Geophys. astropphys. fluid dynamics, 62: 15--36, 1991.



\bibitem{Li}
M. J. Lighthill.
\newblock {\em Studies on magnetohydrodynamic waves and other anisotropic wave motions}.
\newblock Philos. Trans. R. Soc. Lond., Ser. A : 397--430, 1960.

\bibitem{MS}
V. Mazya, and T. Shaposhnikova.
\newblock {\em On the Bourgain, Brezis, and Mironescu theorem concerning limiting embeddings of fractional Sobolev spaces}.
\newblock Journal of Functional Analysis, Vol. 195: 230--238, 2002.






\bibitem{S0}
M. E. Schonbek.
\newblock{\em $L^2$ decay for weak solutions of the Navier-Stokes equations}.
\newblock Archive for Rational Mechanics and Analysis, Vol. 88(3): 209--222, 1985. 

\bibitem{S1}
M. E. Schonbek.
\newblock{\em Large time behavior of solutions of the Navier-Stokes equations}.
\newblock Comm. Partial Differential Equations, Vol. 11: 733--763, 1986. 

\bibitem{S3}
M. E. Schonbek.
\newblock{\em Uniform decay rates for parabolic conservation laws}.
\newblock Journal of Nonlinear Analysis, Vol. 10, No. 9,  943-956, 1986. 

\bibitem{S2}
M. E. Schonbek and M. Wiegner.
\newblock{\em On the decay of higher order norms of the solutions of Navier-Stokes euations}.
\newblock Proc. Roy. Soc. Edinburgh Sect. A 126, 677-685, 1996.


\bibitem{SC}
A. N. Simakov and L. Chacon.
\newblock {\em Quantitative, analytical model for magnetic reconnection in Hall magnetohydrodynamics}.
\newblock Phys. Rev. Lett, 101, 105003, 2008.



\end{thebibliography}
\end{document}